# EKELAND, TAKAHASHI AND CARISTI PRINCIPLES IN PREORDERED QUASI-METRIC SPACES

S. COBZAŞ

ABSTRACT. We prove versions of Ekeland, Takahashi and Caristi principles in preordered quasi-metric spaces, the equivalence between these principles, as well as their equivalence to some completeness results for the underlying quasi-metric space. These extend the results proved in S. Cobzaş, Topology Appl. **265** (2019), 106831, 22, for quasi-metric spaces.

The key tools are Picard sequences for some special set-valued mappings on a preordered quasi-metric space $X$, defined in terms of the preorder and of a function $\varphi$ on $X$.

**Classification MSC 2020:** 46N10 47H10 58E30 54E25 54E35

## 1. INTRODUCTION

Ivar Ekeland announced in 1972, [16] (the proof appeared in 1974, [17]) a theorem asserting the existence of the minimum of a small perturbation of a lower semicontinuous (lsc) function defined on a complete metric space. This result, known as Ekeland Variational Principle (EkVP), proved to be a very versatile tool in various areas of mathematics and applications - optimization theory, geometry of Banach spaces, optimal control theory, economics, social sciences, and others. Some of these applications are presented by Ekeland himself in [18].

At the same time, it turned out that this principle is equivalent to a lot of results in fixed point theory (Caristi fixed point theorem), geometry of Banach spaces (drop property), and others (see [29], for instance). Takahashi [35] (see also [36]) found a sufficient condition for the existence of the minimum of a lsc function on a complete metric space, known as Takahashi minimization principle, which also turned to be equivalent to EkVP (see [35] and [19]).

For convenience, we mention these three principles.

**Theorem 1.1** (Ekeland, Takahashi and Caristi principles). *Let $(X, d)$ be a complete metric space and $\varphi : X \to \mathbb{R} \cup \{\infty\}$ a proper bounded below lsc function. Then the following statements hold.*

(wEk) *There exists $z \in \operatorname{dom} \varphi$ such that $\varphi(z) < \varphi(x) + d(x, z)$ for all $x \in X \smallsetminus \{z\}$.*
(Tak) *If for every $x \in \operatorname{dom} \varphi$ with $\varphi(x) > \inf \varphi(X)$ there exists an element $y \in \operatorname{dom} \varphi \smallsetminus \{x\}$ such that $\varphi(y) + d(x, y) \leq \varphi(x)$, then $\varphi$ attains its minimum on $X$, i.e., there exists $z \in \operatorname{dom} \varphi$ such that $\varphi(z) = \inf \varphi(X)$.*

(Car) *If the mapping $T : X \to X$ satisfies $d(Tx, x) + \varphi(Tx) \leq \varphi(x)$ for all $x \in \operatorname{dom} \varphi$, then $T$ has a fixed point in $\operatorname{dom} \varphi$, i.e., there exists $z \in \operatorname{dom} \varphi$ such that $Tz = z$.*

**Remark 1.2.** Suppose that $\varphi(x) = \infty$. Then, obviously, $\varphi(z) < \varphi(x) + d(y, x)$ for all $z \in \operatorname{dom} \varphi$. In this case, $\varphi(x) > \inf \varphi(X)$ (because $\varphi$ is bounded below). Since the function $\varphi$ is proper, the inequality from (Tak) holds for all $y \in X \smallsetminus \{x\}$. Also the inequality from (Car) is automatically satisfied for such an $x$.

Hence, it suffices to suppose that the conditions from the hypotheses of these results are fulfilled only for points $x \in \operatorname{dom} \varphi$.

The statement (wEk) is called the weak form of the Ekeland variational principle. Later, various versions and extensions of these principles appeared, a good record (up to 2009) being given in the book [24].

Some versions of EkVP and Takahashi minimum principles in $T_1$ quasi-metric spaces were proved in [10] and [2], respectively. In [10] the equivalence of the weak Ekeland variational principle to Caristi fixed point theorem was proved and the implications of the validity of Caristi fixed point theorem principle on the completeness of the underlying quasi-metric space were studied as well. In [2] the same is done for the weak form of Ekeland variational principle and Takahashi minimization principle. The extension of wEk to arbitrary quasi-metric spaces was given in [21] as well as in [14], where the equivalence of these three principles (weak Ekeland, Takahashi and Caristi) is also proved. The case of asymmetric locally convex spaces is treated in [11].

Other extensions, with applications to various areas of social sciences and psychology, were given in [3], [4], [5], [6], [7], [9], [32].

On the other hand, there are many papers dealing with fixed point results for mappings defined on metric spaces endowed with an order, on the line initiated by Ran and Reurings [31] and continued in [25], [26], [27], [30]. These papers contain fixed point results for mappings $f$, defined on a metric space $(X, d)$ endowed with a partial order $\preccurlyeq$, satisfying a weaker contractibility condition, namely: there exists $\alpha \in [0, 1)$ such that

$$(1.1) \qquad \forall x, y \in X, \ x \preccurlyeq y \ \Rightarrow \ d(f(x), f(y)) \leq \alpha \, d(x, y) \,.$$

Some conditions, relating the order and the metric are imposed as, for instance, (2.6). Notice that, if the relation $\preccurlyeq$ is $X \times X$, then (1.1) becomes the usual notion of Banach contraction.

Similar results are obtained for a metric space $(X, d)$ with a graph, see, e.g., Jachymski [20], where a comparison between these two approaches (based on order relations or on graphs) is done. Fixed points on metric spaces endowed with a reflexive binary relation $\mathcal{R}$ were considered in [1] for single-valued maps and in [34] for set-valued mappings. The graphs considered by Jachymski [20], $G = (V(G), E(G))$, where $V(G)$ is the set of vertices and $E(G)$ is the set of edges, satisfy the conditions $V(G) = X$ and $\Delta_X \subseteq E(G) \subseteq X \times X$, so that $E(G)$ is, in fact, a reflexive binary relation on $X$. Here $\Delta_X = \{(x, x) : x \in X\}$ is the diagonal of $X \times X$.

The aim of this paper is to show that results, similar to those proved in [14], hold in a preordered quasi-metric space $(X, d, \preccurlyeq)$, where the preorder $\preccurlyeq$ satisfies (2.6).

In Section 3 we prove versions of Ekeland, Takahashi and Caristi principles in a preordered quasi-metric space as well as their equivalence. Also, the validity of weak Ekeland

Variational Principle in this setting implies a kind of completeness for the underlying space. We conclude this section by showing that in the $T_1$ case one obtains the usual versions of these principles (recall that the topology of a quasi-metric space is only $T_0$).

The key tools used in the proofs are Picard sequences for some special set-valued mappings corresponding to a function $\varphi$ on a preordered quasi-metric space (see Subsection 2.3), which allow a unitary treatment of all these principles. The idea to use Picard sequences appeared in [15] and was subsequently exploited in [8] and [3].

## 2. Quasi-metric spaces

We present some notions and results on quasi-metric spaces needed in the sequel. Full details an other results can be found in [12].

### 2.1. Topological properties.

A *quasi-metric* on an arbitrary set $X$ is a mapping $d : X \times X \to [0, \infty)$ satisfying the following conditions:

$$\begin{aligned}
(\text{QM1}) \quad & d(x, y) \geq 0, \quad \text{and} \quad d(x, x) = 0; \\
(\text{QM2}) \quad & d(x, z) \leq d(x, y) + d(y, z), \\
(\text{QM3}) \quad & d(x, y) = d(y, x) = 0 \Rightarrow x = y,
\end{aligned}$$

for all $x, y, z \in X$. The pair $(X, d)$ is called a *quasi-metric space*. The conjugate of the quasi-metric $d$ is the quasi-metric $\bar{d}(x, y) = d(y, x)$, $x, y \in X$. The mapping $d^s(x, y) = \max\{d(x, y), \bar{d}(x, y)\}$, $x, y \in X$, is a metric on $X$.

If $(X, d)$ is a quasi-metric space, then for $x \in X$ and $r > 0$ we define the balls in $X$ by the formulae

$$\begin{aligned}
B_d(x, r) &= \{y \in X : d(x, y) < r\} \text{ - the open ball, and} \\
B_d[x, r] &= \{y \in X : d(x, y) \leq r\} \text{ - the closed ball.}
\end{aligned}$$

The topology $\tau_d$ (or $\tau(d)$) of a quasi-metric space $(X, d)$ can be defined starting from the family $\mathcal{V}_d(x)$ of neighborhoods of an arbitrary point $x \in X$:

$$\begin{aligned}
V \in \mathcal{V}_d(x) &\iff \exists r > 0 \text{ such that } B_d(x, r) \subseteq V \\
&\iff \exists r' > 0 \text{ such that } B_d[x, r'] \subseteq V.
\end{aligned}$$

The convergence of a sequence $(x_n)$ to $x$ with respect to $\tau_d$, called $d$-convergence and denoted by $x_n \xrightarrow{d} x$, can be characterized in the following way

(2.1) $$x_n \xrightarrow{d} x \iff d(x, x_n) \to 0.$$

Also

(2.2) $$x_n \xrightarrow{\bar{d}} x \iff \bar{d}(x, x_n) \to 0 \iff d(x_n, x) \to 0.$$

As a space equipped with two topologies, $\tau_d$ and $\tau_{\bar{d}}$, a quasi-metric space can be viewed as a bitopological space in the sense of Kelly [22].

The following topological properties are true for quasi-metric spaces.

**Proposition 2.1** (see [12]). *If $(X, d)$ is a quasi-metric space, then the following hold.*
  1. *The ball $B_d(x, r)$ is $\tau_d$-open and the ball $B_d[x, r]$ is $\tau_{\bar{d}}$-closed. The ball $B_d[x, r]$ need not be $\tau_d$-closed.*

2. *The topology $\tau_d$ is $T_0$. It is $T_1$ if and only if $d(x,y) > 0$ for all $x \neq y$ in $X$. The topology $\tau_d$ is $T_2$ (or Hausdorff) if and only if every convergent sequence has a unique limit.*
3. *For every fixed $x \in X$, the mapping $d(x, \cdot) : X \to (\mathbb{R}, |\cdot|)$ is $\tau_d$-usc and $\tau_{\bar{d}}$-lsc. For every fixed $y \in X$, the mapping $d(\cdot, y) : X \to (\mathbb{R}, |\cdot|)$ is $\tau_d$-lsc and $\tau_{\bar{d}}$-usc.*

Notice that

$$y \in \overline{\{x\}} \iff d(y,x) = 0. \tag{2.3}$$

Indeed,

$$\begin{aligned} y \in \overline{\{x\}} &\iff \forall r > 0,\, x \in B(y,r) \\ &\iff \forall r > 0,\, d(y,x) < r \\ &\iff d(y,x) = 0. \end{aligned}$$

The lack of symmetry in the definition of quasi-metric spaces causes a lot of troubles, mainly concerning completeness, compactness and total boundedness in such spaces. There are a lot of completeness notions in quasi-metric spaces, all agreeing with the usual notion of completeness in the metric case, each of them having its advantages and weaknesses (see [33], or [12]).

As in what follows we shall work only with two of these notions, we present only them, referring to [12] for others.

We use the notation

$$\mathbb{N} = \{1, 2, \dots\} \text{ -- the set of natural numbers,}$$
$$\mathbb{N}_0 = \mathbb{N} \cup \{0\} \text{ -- the set of non-negative integers.}$$

**Definition 2.2.** Let $(X, d)$ be a quasi-metric space. A sequence $(x_n)$ in $(X, d)$ is called:
  (a) *left $d$-$K$-Cauchy* if for every $\varepsilon > 0$ there exists $n_\varepsilon \in \mathbb{N}$ such that

$$\begin{aligned} &\forall n, m, \text{ with } n_\varepsilon \leq n < m, \quad d(x_n, x_m) < \varepsilon \\ &\iff \forall n \geq n_\varepsilon,\, \forall k \in \mathbb{N}, \quad d(x_n, x_{n+k}) < \varepsilon; \end{aligned} \tag{2.4}$$

  (b) *right $d$-$K$-Cauchy* if for every $\varepsilon > 0$ there exists $n_\varepsilon \in \mathbb{N}$ such that

$$\begin{aligned} &\forall n, m, \text{ with } n_\varepsilon \leq n < m, \quad d(x_m, x_n) < \varepsilon \\ &\iff \forall n \geq n_\varepsilon,\, \forall k \in \mathbb{N}, \quad d(x_{n+k}, x_n) < \varepsilon. \end{aligned} \tag{2.5}$$

The quasi-metric space $(X, d)$ is called:
- *sequentially left $d$-$K$-complete* if every left $d$-$K$-Cauchy is $d$-convergent;
- *sequentially right $d$-$K$-complete* if every right $d$-$K$-Cauchy is $d$-convergent.

**Remarks 2.3.**
  1. It is obvious that a sequence is left $d$-$K$-Cauchy if and only if it is right $\bar{d}$-$K$-Cauchy.
  2. Some examples of $d$-convergent sequences which are not left $d$-$K$-Cauchy, show that in the asymmetric case the situation is far more complicated than in the symmetric one (see [33]).

3. If each convergent sequence in a regular quasi-metric space $(X, d)$ admits a left $K$-Cauchy subsequence, then $X$ is metrizable (see [23]).

We mention the following result on Cauchy sequences.

**Proposition 2.4.** *Let $(X, d)$ be a quasi-metric space. If a right $K$-Cauchy sequence $(x_n)$ contains a subsequence convergent to some $x \in X$, then the sequence $(x_n)$ converges to $x$.*

**Convention.** *In the following, when speaking about metric or topological properties in a quasi-metric space $(X, d)$ we shall always understand those corresponding to $d$ and we shall omit $d$ or $\tau_d$, i.e., we shall write "$(x_n)$ is right $K$-Cauchy" instead of "$(x_n)$ is right $d$-$K$-Cauchy", $\overline{A}$ instead of $\overline{A}^d$, etc.*

2.2. **Preordered quasi-metric spaces.** We consider a preorder $\preccurlyeq$ on a quasi-metric space $(X, d)$ (i.e. a reflexive and transitive relation). If $\preccurlyeq$ is also reflexive, then it is called and order. Consider the following property relating the preorder $\preccurlyeq$ and the quasi-metric $d$: for every sequence $(x_n)_{n \in \mathbb{N}}$ in $X$

(2.6) $\quad x_n \preccurlyeq x_{n+1}$, for all $n \in \mathbb{N}$, and $x_n \xrightarrow{d} x \implies x_n \preccurlyeq x$ for all $n \in \mathbb{N}$.

We call a sequence $(x_n)_{n \in \mathbb{N}}$ satisfying the condition $x_n \preccurlyeq x_{n+1}$, $n \in \mathbb{N}$, *increasing*.

**Remark 2.5.** This condition was introduced in [26] in the study of fixed points for contractions in ordered metric spaces. Some authors consider a weakened version of (2.6), where one asks that the conclusion holds only for a subsequence $(x_{n_k})_{k \in \mathbb{N}}$ of $(x_n)$. As remarked Jachymski [20], the transitivity of $\preccurlyeq$ implies the equivalence of these conditions.

Also, if $\preccurlyeq$ is a reflexive relation on $X$ satisfying (2.6), then it is transitive.

Indeed, if $x_n \preccurlyeq x_{n+1}$ and $x_{n_k} \leq x$, $k \in \mathbb{N}$, for some subsequence $(x_{n_k})_{k \in \mathbb{N}}$ of $(x_n)$, then $k \leq n_k$ implies $x_k \preccurlyeq x_{n_k} \preccurlyeq x$, for all $k \in \mathbb{N}$.

Suppose that $\preccurlyeq$ is a reflexive relation on $X$ satisfying (2.6) and let $x \preccurlyeq y$ and $y \preccurlyeq z$. Considering the sequence $x_1 = x$, $x_2 = y$ and $x_n = z$ for $n \geq 3$, it follows that $(x_n)_{n \in \mathbb{N}}$ is an increasing sequence convergent to $z$, so that $x = x_1 \preccurlyeq z$.

We introduce now some notions expressed in terms of increasing sequences.

**Definition 2.6.** Let $(X, d, \preccurlyeq)$ be a preordered quasi-metric space.
- A sequence $(x_n)_{n \in \mathbb{N}}$ in $X$ is called *increasingly left-K-Cauchy* (*right-K-Cauchy*) if it is increasing and satisfies the condition (a) (resp. (b)) from Definition 2.2.
- The space $X$ is called *increasingly left* (*right*) K-*complete* if every increasingly left (right) K-Cauchy sequence is convergent.
- A subset $Y$ of $X$ is called *increasingly closed* if it contains the limit of every increasing convergent sequence $(y_n)$ in $Y$.
- A function $\varphi : X \to \mathbb{R} \cup \{\infty\}$ is called *increasingly lower semi-continuous* (lsc) if

(2.7) $$\varphi(x) \leq \liminf_n \varphi(x_n),$$

for every increasing sequence $(x_n)_{n \in \mathbb{N}}$ in $X$ convergent to $x$.

**Remark 2.7.** In [3, Def. 9] a different notion of decreasingly lsc function was considered. Namely, a function $\varphi : X \to \mathbb{R} \cup \{\infty\}$ is called decreasingly lsc if $\varphi$ satisfies (2.7) for every sequence $(x_n)$ in $X$ such that $(x_n)$ is $\bar{d}$-convergent to $x$ and $\varphi(x_{n+1}) \leq \varphi(x_n)$, $n \in \mathbb{N}$. Obviously, an "increasingly" version of this notion can be defined for the $d$-convergent sequences. Another version of lower semicontinuity, called nearly lower semicontinuity, was defined in [21] supposing that (2.7) holds only for sequences with pairwise distinct terms. A characterization of this property in terms of the monotonicity with respect to the specialization order ($x \leq_d y \iff d(x,y) = 0$) was given in [14].

2.3. **Picard sequences in preordered quasi-metric spaces.** Let $(X, d, \preccurlyeq)$ be a preordered quasi-metric space and $\varphi : X \to \mathbb{R} \cup \{\infty\}$ be a function. We introduce another preorder $\leq_\varphi$ on $X$ by

$$(2.8) \qquad x \leq_\varphi y \iff \varphi(y) + d(y, x) \leq \varphi(x),$$

for $x, y \in X$.

It is easy to check that $\leq_\varphi$ is a preorder on $X$ (i.e. reflexive and transitive) and an order on $\operatorname{dom}\varphi := \{x \in X : \varphi(x) < \infty\}$.

The reflexivity is obvious.

If $\varphi(x) = \infty$, then $x \leq_\varphi y$ for all $y \in X$, so that the transitivity holds trivially.

If $x \in \operatorname{dom}\varphi$ and $x \leq_\varphi y$, then $y \in \operatorname{dom}\varphi$. Hence

$$\begin{aligned} x \leq_\varphi y &\iff d(y, x) \leq \varphi(x) - \varphi(y) \text{ and} \\ y \leq_\varphi z &\iff d(z, y) \leq \varphi(y) - \varphi(z), \end{aligned}$$

so that

$$d(z, x) \leq d(z, y) + d(y, x) \leq \varphi(x) - \varphi(z),$$

that is $x \leq_\varphi z$.

Also,

$$\begin{aligned} x \leq_\varphi y &\iff d(y, x) \leq \varphi(x) - \varphi(y) \text{ and} \\ y \leq_\varphi x &\iff d(x, y) \leq \varphi(y) - \varphi(x), \end{aligned}$$

so that

$$0 \leq d(y, x) + d(x, y) \leq 0.$$

Hence, $d(y, x) = d(x, y) = 0$ and so, by (QM3), $x = y$.

For $x \in \operatorname{dom}\varphi$ consider the set $S(x) \subseteq X$ given by

$$(2.9) \qquad S(x) = \{y \in X : x \preccurlyeq y \text{ and } x \leq_\varphi y\}.$$

**Remark 2.8.** For $x \in X \smallsetminus \operatorname{dom}\varphi$, $S(x) = \{y \in X : x \preccurlyeq y\}$. It is worth to mention that sets of this kind were used by Penot [28] as early as 1977 in a proof of Caristi fixed point theorem in complete metric spaces.

**Proposition 2.9.** *Let $x \in \operatorname{dom} \varphi$. The following properties hold.*

  (i) $x \in S(x)$   *and*   $S(x) \subseteq \operatorname{dom} \varphi$;
  (ii) $y \in S(x) \Rightarrow \varphi(y) \leq \varphi(x)$ *and* $S(y) \subseteq S(x)$;
  (iii) $y \in S(x) \smallsetminus \overline{\{x\}} \Rightarrow \varphi(y) < \varphi(x)$;
       $y \in S(x)$ *and* $\varphi(y) = \varphi(x) \Rightarrow y \in \overline{\{x\}}$;
  (iv) *if $\varphi$ is bounded below, then*
$$S(x) \smallsetminus \overline{\{x\}} \neq \emptyset \Rightarrow \varphi(x) > \inf \varphi(S(x));$$
  (v) *if the preorder $\preccurlyeq$ satisfies (2.6) and $\varphi$ is increasingly lsc, then $S(x)$ is increasingly closed, i.e., $y \in S(x)$ for every increasing sequence $(y_n)$ in $S(x)$ convergent to $y \in X$.*

*Proof.* The relations (i) are immediate consequences of the definition of $S(x)$.

(ii) If $y \in S(x)$, then $0 \leq d(y,x) \leq \varphi(x) - \varphi(y)$ implies $\varphi(y) \leq \varphi(x)$. The inclusion $S(y) \subseteq S(x)$ follows by the transitivity of the relations $\preccurlyeq$ and $\leq_\varphi$.

(iii)  By (2.3), $d(y,x) > 0$ for every $y \in S(x) \smallsetminus \overline{\{x\}}$, so that $x \leq_\varphi y$ yields
$$0 < d(y,x) \leq \varphi(x) - \varphi(y),$$
that is, $\varphi(y) < \varphi(x)$.

The second implication from (iii) follows from the first one.

(iv) If there exists $y \in S(x) \smallsetminus \overline{\{x\}}$, then, by (iii),
$$\inf \varphi(S(x)) \leq \varphi(y) < \varphi(x).$$

(v) If $y_n \in S(x)$, $y_n \preccurlyeq y_{n+1}$, $n \in \mathbb{N}$, and $y_n \xrightarrow{d} y$, then, by (2.6), $x \preccurlyeq y_1 \preccurlyeq y$. The function $d(\cdot,x)$ is $d$-lsc (see Proposition 2.1.3) and $\varphi$ is increasingly $d$-lsc, so that
$$\varphi(y_n) + d(y_n,x) \leq \varphi(x), \ n \in \mathbb{N},$$
implies
$$\varphi(y) + d(y,x) \leq \liminf_n [\varphi(y_n) + d(y_n,x)] \leq \varphi(x),$$
that is $x \leq_\varphi y$. Consequently $y \in S(x)$. $\square$

**Remark 2.10.** If $\inf \varphi(S(x)) = -\infty$, then $\varphi(x) > \inf \varphi(S(x))$ no matter that $S(x) \smallsetminus \overline{\{x\}}$ is nonempty or not.

A *Picard sequence* corresponding to a set-valued mapping $F : X \rightrightarrows X$ is a sequence $(x_n)_{n \in \mathbb{N}_0}$ such that
$$x_{n+1} \in F(x_n) \ \text{ for all } \ n \in \mathbb{N}_0,$$
for a given initial point $x_0 \in X$. This notion was introduced in [15] (see also [8]). We shall work with some special Picard sequences corresponding to the set-valued mapping $S : \operatorname{dom} \varphi \rightrightarrows \operatorname{dom} \varphi$, where $S(x)$ is given by (2.9) for all $x \in \operatorname{dom} \varphi$. Put
$$J(x) = \inf \varphi(S(x)).$$

Observe that if, for some $x \in \mathrm{dom}\,\varphi$, $J(x) \in \mathbb{R}$ and $\varphi(x) > J(x)$, then there exists an element $y \in S(x)$ such that
$$\varphi(y) < \frac{\varphi(x) + J(x)}{2}\,.$$

**Proposition 2.11** (Picard sequences). *Let $(X, d, \preccurlyeq)$ be a preordered quasi-metric space and $\varphi : X \to \mathbb{R} \cup \{\infty\}$ be a bounded below proper function.*

*Let $x_0 \in \mathrm{dom}\,\varphi$. We distinguish two situations.*

1. *There exists $m \in \mathbb{N}_0$ such that*

(2.10)
  (i) $\varphi(x_k) > J(x_k)$,
  (ii) $x_{k+1} \in S(x_k)$ and $\varphi(x_{k+1}) < (\varphi(x_k) + J(x_k))/2$,
      for all $0 \leq k \leq m - 1$, and
  (iii) $\varphi(x_m) = J(x_m)$.

*Putting $z = x_m$, the following conditions are satisfied:*

(2.11)
  (i) $x_k \preccurlyeq x_{k+1}$, $S(x_{k+1}) \subseteq S(x_k)$ and $\varphi(x_{k+1}) < \varphi(x_k)$,
      for all $0 \leq k \leq m - 1$;
  (ii) $z \in S(x_k)$ and $S(z) \subseteq S(x_k)$ for $0 \leq k \leq m$;
  (iii) $\varphi(y) = \varphi(z) = J(z)$ and $S(y) \subseteq \overline{\{y\}}$,
      for all $y \in S(z)$.

2. *There exists a sequence $(x_n)_{n \in \mathbb{N}_0}$ such that*

(2.12)
  (i) $\varphi(x_n) > J(x_n)$;
  (ii) $x_{n+1} \in S(x_n)$ and $\varphi(x_{n+1}) < (\varphi(x_n) + J(x_n))/2$,

*for all $n \in \mathbb{N}_0$.*

*Then the sequence $(x_n)_{n \in \mathbb{N}_0}$ satisfies the conditions*

(2.13)
  (i) $x_n \preccurlyeq x_{n+1}$, $S(x_{n+1}) \subseteq S(x_n)$ and $\varphi(x_{n+1}) < \varphi(x_n)$,
      for all $n \in \mathbb{N}_0$;
  (ii) there exist the limits $\alpha := \lim_{n \to \infty} \varphi(x_n) = \lim_{n \to \infty} J(x_n) \in \mathbb{R}$;
  (iii) $x_{n+k} \in S(x_n)$ for all $n, k \in \mathbb{N}_0$;
  (iv) $(x_n)_{n \in \mathbb{N}_0}$ is increasingly right $K$-Cauchy.

*If, in addition, the space $X$ is sequentially increasingly right $K$-complete, the preorder $\preccurlyeq$ satisfies the condition (2.6) and the function $\varphi$ is increasingly lsc, then the sequence $(x_n)_{n \in \mathbb{N}_0}$ is convergent to a point $z \in X$ such that*

(2.14)
  (i) $z \in S(x_n)$ and $S(z) \subseteq S(x_n)$ for all $n \in \mathbb{N}_0$;
  (ii) $\varphi(y) = \varphi(z) = J(z) = \alpha$ and $S(y) \subseteq \overline{\{y\}}$, for all $y \in S(z)$,
      where $\alpha$ is given by (2.13).(ii).

*Proof.* Suppose that we have found $x_0, x_1, \ldots, x_m$ satisfying the conditions (i) and (ii) from (2.10). If $\varphi(x_m) = J(x_m)$ then $x_0, x_1, \ldots, x_m$ satisfy (2.10).

If $\varphi(x_m) > J(x_m)$, then there exists $x_{m+1} \in S(x_m)$ such that $\varphi(x_{m+1}) < (\varphi(x_m) + J(x_m))/2$. Supposing that this procedure continues indefinitely, we find a sequence $(x_n)_{n \in \mathbb{N}_0}$ satisfying (2.12).

1. Suppose that $x_0, x_1, \ldots, x_m$ satisfy (2.10) and let $z = x_m$. Then, by Proposition 2.9.(ii), $x_{k+1} \in S(x_k)$ implies $S(x_{k+1}) \subseteq S(x_k)$ for $k = 0, 1, \ldots, m-1$.

If $y \in S(z)$, then, by Proposition 2.9.(ii),
$$J(z) \leq \varphi(y) \leq \varphi(z) = J(z).$$

It follows $\varphi(y) = \varphi(z) = J(z)$ for all $y \in S(z)$.

If $y \in S(z)$ and $x \in S(y) \subseteq S(z)$, then $\varphi(x) = \varphi(z) = \varphi(y)$, so that, by Proposition 2.9.(iii), $x \in \overline{\{y\}}$.

We have shown that all of the conditions from (2.11) are satisfied by $z$.

2. Suppose now that the sequence $(x_n)_{n \in \mathbb{N}_0}$ satisfies (2.12).

By Proposition 2.9.(ii), the relation $x_{n+1} \in S(x_n)$ implies
$$x_n \preccurlyeq x_{n+1} \quad \text{and} \quad S(x_{n+1}) \subseteq S(x_n),$$
for all $n \in \mathbb{N}_0$.

Also, by (2.12).(ii),
$$\varphi(x_{n+1}) < \varphi(x_n) \quad \text{for all} \ \ n \in \mathbb{N}_0,$$
so that (2.13).(i) holds.

Since $\varphi$ is bounded below and, by (2.13).(i), $(\varphi(x_n))_{n \in \mathbb{N}_0}$ is strictly decreasing, there exists the limit
$$\alpha := \lim_{n \to \infty} \varphi(x_n) = \inf_{n \in \mathbb{N}_0} \varphi(x_n) \in \mathbb{R}.$$

By (2.12),
$$2\varphi(x_{n+1}) - \varphi(x_n) < J(x_n) < \varphi(x_n),$$
for all $n \in \mathbb{N}_0$. Letting $n \to \infty$, one obtains

(2.15) $$\lim_{n \to \infty} J(x_n) = \lim_{n \to \infty} \varphi(x_n) = \alpha.$$

By (2.13).(i),
$$x_{n+k} \in S(x_{n+k}) \subseteq S(x_n),$$
proving (2.13).(iii). Hence $x_n \leq_\varphi x_{n+k}$, that is,

(2.16) $$d(x_{n+k}, x_n) \leq \varphi(x_n) - \varphi(x_{n+k}) \quad \text{for all} \quad n, k \in \mathbb{N}_0.$$

Since the sequence $(\varphi(x_n))_{n \in \mathbb{N}_0}$ is Cauchy, this implies that $(x_n)_{n \in \mathbb{N}_0}$ is right $K$-Cauchy. It is also increasing.

Suppose now that $X$ is sequentially increasingly right $K$-complete and $\varphi$ is increasingly lsc. Then there exists $z \in X$ such that $x_n \xrightarrow{d} z$.

Since the sequence $(x_n)_{n \in \mathbb{N}_0}$ is increasing and $d$-convergent to $z$, condition (2.6) implies
$$x_n \preccurlyeq z \quad \text{for all} \quad n \in \mathbb{N}_0.$$

Taking into account the facts that the function $d(\cdot, x_n)$ is lsc (Proposition 2.1.3) and $\varphi$ is increasingly lsc, the inequalities (2.16) yield

$$\varphi(z) + d(z, x_n) \leq \liminf_{k \to \infty} [\varphi(x_{n+k}) + d(x_{n+k}, x_n)] \leq \varphi(x_n).$$

Consequently, $z \in S(x_n)$ and so $S(z) \subseteq S(x_n)$, for all $n \in \mathbb{N}_0$.

If $y \in S(z) \subseteq S(x_n)$, then

$$J(x_n) \leq J(z) \leq \varphi(y) \leq \varphi(y) + d(y, x_n) \leq \varphi(x_n),$$

that is

$$J(x_n) \leq J(z) \leq \varphi(y) \leq \varphi(x_n),$$

for all $n \in \mathbb{N}_0$.

Letting $n \to \infty$ and taking into account (2.13).(ii), one obtains

$$\varphi(y) = J(z),$$

for all $y \in S(z)$.

The inclusion from (2.14).(ii) follows from Proposition 2.9.(iii), as in the case of the corresponding one from (2.11).(iii). □

**Remark 2.12.** If the function $\varphi$ satisfies

$$\varphi(x) > \inf \varphi(S(x)) \quad \text{for all} \quad x \in \operatorname{dom} \varphi,$$

then, for every $x_0 \in \operatorname{dom} \varphi$, there exists a Picard sequence $(x_n)_{n \in \mathbb{N}_0}$ satisfying (2.12) (and so (2.13), as well).

3. EKELAND, TAKAHASHI AND CARISTI PRINCIPLES IN PREORDERED QUASI-METRIC SPACES

Ekeland, Takahashi and Caristi principles in metric spaces (see Theorem 1.1) can be expressed in terms of some sets, similar to $S(x)$ (see (2.9)), in the following form.

**Theorem 3.1.** *Let $(X, d)$ be a complete metric space, $\varphi : X \to \mathbb{R} \cup \{\infty\}$ be a proper bounded below lsc function and let*

$$R(x) = \{y \in X : \varphi(y) + d(y, x) \leq \varphi(x)\},$$

*for $x \in \operatorname{dom} \varphi$. Then the following results hold.*

(wEk) *There exists $z \in \operatorname{dom} \varphi$ such that $R(z) = \{z\}$.*

(Tak) *If*

$$\varphi(x) > \inf \varphi(X) \implies R(x) \smallsetminus \{x\} \neq \emptyset,$$

*for all $x \in \operatorname{dom} \varphi$, then $\varphi$ attains its minimum on $X$, i.e., there exists $z \in X$ such that $\varphi(z) = \inf \varphi(X)$.*

(Car) *If the mapping $T : X \to X$ satisfies $Tx \in R(x)$ for all $x \in \operatorname{dom} \varphi$, then $T$ has a fixed point in $X$, i.e., there exists $z \in X$ such that $Tz = z$.*

**Remark 3.2.** If $\varphi(x) = \infty$, then $R(x) = X \neq \{x\}$, because $\varphi$ is proper. Also the condition $Tx \in R(x)$ automatically holds, so it suffices to ask that the corresponding conditions are fulfilled only for $x \in \operatorname{dom} \varphi$.

In [14] these results were extended to quasi-metric spaces. In this paper we go further and prove some versions of these results in preordered quasi-metric spaces. Along this section we shall use the notation: $(X, d, \preccurlyeq)$ will be a preordered quasi-metric space, $\varphi : X \to \mathbb{R} \cup \{\infty\}$ a function and $\leq_\varphi$ the order defined by (2.8). Also, for $x \in X$ put

(3.1) $\qquad S(x) = \{y \in X : x \preccurlyeq y \text{ and } x \leq_\varphi y\} \quad \text{and} \quad J(x) = \inf \varphi(S(x)).$

## 3.1. Ekeland variational principle.
We start with a version of weak Ekeland principle.

**Theorem 3.3.** *Let $(X, d, \preccurlyeq)$ be a sequentially increasingly right $K$-complete preordered quasi-metric space and let $\varphi : X \to \mathbb{R} \cup \{\infty\}$ be a proper bounded below increasingly lsc function. Suppose that the preorder $\preccurlyeq$ satisfies (2.6). Then there exists $z \in \operatorname{dom} \varphi$ such that*

(3.2) $\qquad \varphi(y) = \varphi(z) \ \text{ for all } \ y \in S(z).$

*In this case it follows that, for every $y \in S(z)$,*

(3.3) (i) $S(y) \subseteq \overline{\{y\}}$ *and*
(ii) $\varphi(y) < \varphi(x) + d(x, y), \forall x \in \operatorname{dom} \varphi \smallsetminus S(y)$ *with* $y \preccurlyeq x$ *and* $\forall x \in X \smallsetminus \operatorname{dom} \varphi$.

*Proof.* By (2.11).(ii) and (2.14).(ii), there exists $z \in \operatorname{dom} \varphi$ satisfying (3.2) and (3.3).(i).

Recall that $S(z) \subseteq \operatorname{dom} \varphi$ and $S(y) \subseteq S(z)$. If $x \in \operatorname{dom} \varphi \smallsetminus S(y)$ and $y \preccurlyeq x$, then, by the definition of the set $S(y)$,

$$\varphi(y) < \varphi(x) + d(x, y),$$

so that (3.3).(ii) holds for this $x$. If $x \in X \smallsetminus \operatorname{dom} \varphi$, then $\varphi(x) = \infty$, so that (3.3).(ii) holds for this $x$ too. $\square$

The full version of Ekeland variational principle is the following.

**Theorem 3.4.** *Let $(X, d, \preccurlyeq)$ be a preordered quasi-metric space. Suppose that $X$ is sequentially increasingly right $K$-complete and that $\preccurlyeq$ satisfies (2.6). Let $\varphi : X \to \mathbb{R} \cup \{\infty\}$ be a proper, bounded below and increasingly lsc function. If $\varepsilon, \lambda > 0$ and $x_0 \in X$ is such that*

(3.4) $\qquad \varphi(x_0) \leq \varepsilon + \inf \varphi(X),$

*then there exists $z \in \operatorname{dom} \varphi$ such that*

(3.5)
(i) $\varphi(z) + \dfrac{\varepsilon}{\lambda} d(z, x_0) \leq \varphi(x_0) \quad (\text{and so } \varphi(z) \leq \varphi(x_0));$
(ii) $d(z, x_0) \leq \lambda;$
(iii) $\varphi(y) = \varphi(z) \ \text{ for all } \ y \in S(z);$
(iv) $\varphi(z) < \varphi(x) + \dfrac{\varepsilon}{\lambda} d(x, z),$
$\forall x \in \operatorname{dom} \varphi \smallsetminus S(z) \ \text{ with } \ z \preccurlyeq x \text{ and } \forall x \in X \smallsetminus \operatorname{dom} \varphi.$

*Proof.* For convenience, put $\gamma = \varepsilon/\lambda$ and $d_\gamma = \gamma d$. Then $d_\gamma$ is a quasi-metric on $X$ Lipschitz equivalent to $d$, so that $(X, d_\gamma)$ is also sequentially increasingly right $K$-complete.

Let

$$X_0 = \{x \in X : x_0 \preccurlyeq x \ \text{ and } \ \varphi(x) \leq \varphi(x_0) + d_\gamma(x_0, x)\}.$$

Remark that $X_0 \subseteq \operatorname{dom} \varphi$.

*Claim* I. *The set $X_0$ is increasingly closed and $x_0 \in X_0$.*

Indeed, let $(x_n)$ be an increasing sequence in $X_0$ which is $d_\gamma$-convergent to some $x \in X$, i.e., $\lim_{n\to\infty} d_\gamma(x, x_n) = 0$. Then, by (2.6),

$$x_0 \preccurlyeq x_1 \preccurlyeq x,$$

that is $x_0 \preccurlyeq x$.

Also, $x_n \in X_0$ implies

$$\begin{aligned}\varphi(x_n) &\leq \varphi(x_0) + d_\gamma(x_0, x_n) \\ &\leq \varphi(x_0) + d_\gamma(x_0, x) + d_\gamma(x, x_n),\end{aligned}$$

for all $n \in \mathbb{N}$. Taking into account that the function $\varphi$ is increasingly lsc, one obtains

$$\varphi(x) \leq \liminf_{n\to\infty} \varphi(x_n) \leq \varphi(x_0) + d_\gamma(x_0, x).$$

Consequently, $x \in X_0$. It is obvious that $x_0 \in X_0$.

For $y \in X_0$ put

$$S_{X_0}(y) := \{x \in X_0 : y \preccurlyeq x \text{ and } \varphi(x) + d_\gamma(x, y) \leq \varphi(y)\} = X_0 \cap S(y),$$

where

$$S(y) := \{z \in X : y \preccurlyeq z \text{ and } \varphi(z) + d_\gamma(z, y) \leq \varphi(y)\}.$$

*Claim* II. *For every $y \in X_0$*

(3.6) $\quad$ (i) $\quad S(y) \subseteq X_0$ *so that* $S_{X_0}(y) = S(y),$ *and*

$\qquad$ (ii) $\quad \varphi(y) < \varphi(x) + d_\gamma(x, y)$ *for every* $x \in \operatorname{dom} \varphi \smallsetminus X_0$ *with* $x_0 \preccurlyeq x$.

We have

$$S(y) \subseteq X_0 \iff S(y) \smallsetminus X_0 = \emptyset.$$

Suppose that $S(y) \smallsetminus X_0 \neq \emptyset$. Observe that $y \in X_0$ and $z \in S(y)$ imply $x_0 \preccurlyeq y$ and $y \preccurlyeq z$, so that $x_0 \preccurlyeq z$. Hence, if $z \in S(y) \smallsetminus X_0$, then the following inequalities are satisfied:

$$\begin{aligned}\varphi(y) &\leq \varphi(x_0) + d_\gamma(x_0, y) & (y \in X_0), \\ \varphi(z) + d_\gamma(z, y) &\leq \varphi(y) & (z \in S(y)), \\ \varphi(x_0) + d_\gamma(x_0, z) &< \varphi(z) & (z \notin X_0).\end{aligned}$$

These lead to the contradiction:

$$\begin{aligned}\varphi(z) + d_\gamma(z, y) &\leq \varphi(y) \leq \varphi(x_0) + d_\gamma(x_0, y) \\ &\leq \varphi(x_0) + d_\gamma(x_0, z) + d(z, y) \\ &< \varphi(z) + d_\gamma(z, y).\end{aligned}$$

If $x \in \operatorname{dom} \varphi \smallsetminus X_0$ and $x_0 \preccurlyeq x$, then

$$\varphi(x_0) + d_\gamma(x_0, x) < \varphi(x),$$

so that
$$\varphi(y) \leq \varphi(x_0) + d_\gamma(x_0, y)$$
$$\leq \varphi(x_0) + d_\gamma(x_0, x) + d_\gamma(x, y)$$
$$< \varphi(x) + d_\gamma(x, y),$$

i.e, (3.6).(ii) holds too.

By Claim I, $(X_0, d_\gamma, \preccurlyeq)$ is sequentially increasingly right $K$-complete. Applying Proposition 2.11 to this space we find an element $z \in X_0$ such that

(3.7) $$\varphi(y) = \varphi(z) \text{ for all } y \in S_{X_0}(z) = S(z).$$

This shows that $z$ satisfies (3.5).(iii).

Now, by (2.11).(ii) and (2.14).(i), $z \in S_{X_0}(x_0) = S(x_0)$, which is equivalent to (3.5).(i).

Also, (3.5).(i) and (3.4) imply

$$\frac{\varepsilon}{\lambda} d(z, x_0) \leq \varphi(x_0) - \varphi(z) \leq \varphi(x_0) - \inf \varphi(X) \leq \varepsilon,$$

so that
$$d(z, x_0) \leq \lambda,$$

i.e., (3.5).(ii) holds too.

The inequality (3.5).(iv) holds for arbitrary $x \in X \smallsetminus \operatorname{dom} \varphi$ and, by the definition of the set $S(z)$, for every $x \in \operatorname{dom} \varphi \smallsetminus S(z)$ with $z \preccurlyeq x$. □

## 3.2. Takahashi principle.
The following form of Takahashi principle holds.

**Theorem 3.5.** *Let $(X, d, \preccurlyeq)$ be a sequentially increasingly right $K$-complete preordered quasi-metric space and let $\varphi : X \to \mathbb{R} \cup \{\infty\}$ be a proper bounded below increasingly lsc function. Suppose that the preorder $\preccurlyeq$ satisfies (2.6). If the following condition*

(3.8) $$\varphi(x) > \inf \varphi(X) \implies \exists y \in S(x), \varphi(y) < \varphi(x),$$

*holds for every $x \in \operatorname{dom} \varphi$, then the function $\varphi$ attains its minimum on $X$, i.e, there exists $z \in X$ such that $\varphi(z) = \inf \varphi(X)$.*

*Proof.* Suppose, by contradiction, that

(3.9) $$\varphi(x) > \inf \varphi(X),$$

for all $x \in \operatorname{dom} \varphi$. Then, by (3.8),

(3.10) $$\forall x \in \operatorname{dom} \varphi, \exists y \in S(x), \varphi(y) < \varphi(x),$$

or, equivalently,

(3.11) $$\varphi(x) > \inf \varphi(S(x)),$$

for all $x \in \operatorname{dom} \varphi$.

Let $x_0 \in \operatorname{dom} \varphi$. By (3.11), Remark 2.12 and Proposition 2.11.2, there exists $z \in X$ such that (2.14) holds. By (3.10) there exists $y \in S(z)$ with

(3.12) $$\varphi(y) < \varphi(z),$$

in contradiction to (2.14).(ii). □

**Corollary 3.6.** *Suppose that $(X, d, \preccurlyeq)$ and $\varphi : X \to \mathbb{R} \cup \{\infty\}$ satisfy the hypotheses of Theorem 3.5. If, for every $x \in \operatorname{dom} \varphi$,*

$$\tag{3.13} \varphi(x) > \inf \varphi(X) \;\Rightarrow\; S(x) \smallsetminus \overline{\{x\}} \neq \emptyset ,$$

*then the function $\varphi$ attains its minimum on $X$.*

*Proof.* Condition (3.13) means that, for every $x \in X$,

$$\varphi(x) > \inf \varphi(X) \;\Rightarrow\; \exists y \in S(x) \smallsetminus \overline{\{x\}} .$$

By Proposition 2.9.(iii), this implies

$$\varphi(y) < \varphi(x) ,$$

so we can apply Theorem 3.5 to conclude. $\square$

### 3.3. Caristi fixed point theorem.
In the $T_0$ setting we can prove only some weaker forms of Caristi fixed point theorem for single-valued and set-valued mappings. If $X$ is $T_1$, then the proper versions hold (see Subsection 3.5).

**Theorem 3.7** (Caristi's theorem). *Suppose that $(X, d, \preccurlyeq)$ is a sequentially increasingly right $K$-complete preordered quasi-metric space such that the preorder $\preccurlyeq$ satisfies (2.6). Let $\varphi : X \to \mathbb{R} \cup \{\infty\}$ be a proper, bounded below and increasingly lsc function.*

*1. If the mapping $T : X \to X$ satisfies*

$$\tag{3.14} x \preccurlyeq Tx \quad \text{and} \quad d(Tx, x) + \varphi(Tx) \le \varphi(x) ,$$

*for all $x \in \operatorname{dom} \varphi$, then there exists $z \in \operatorname{dom} \varphi$ such that*

$$\tag{3.15} \varphi(Tz) = \varphi(z) \quad \text{and} \quad Tz \in \overline{\{z\}}.$$

*2. If $T : X \rightrightarrows X$ is a set-valued mapping such that*

$$\tag{3.16} S(x) \cap Tx \neq \emptyset ,$$

*for every $x \in \operatorname{dom} \varphi$, then there exists $z \in \operatorname{dom} \varphi$ such that*

$$\tag{3.17} \varphi(z) \in \varphi(Tz) \quad \text{and} \quad Tz \cap \overline{\{z\}} \neq \emptyset.$$

*Proof.* Observe that condition (3.14) is equivalent to

$$\tag{3.18} Tx \in S(x) ,$$

for all $x \in \operatorname{dom} \varphi$. This shows that (3.16) becomes (3.14) for single-valued mappings, so it suffices to prove only the statement 2.

By (2.11)(iii) and (2.14)(ii), there exists $z \in \operatorname{dom} \varphi$ with $S(z) \subseteq \overline{\{z\}}$ and $\varphi(y) = \varphi(z)$ for all $y \in S(z)$. By (3.16),

$$Tz \cap \overline{\{z\}} \supseteq Tz \cap S(z) \neq \emptyset .$$

Also,

$$\varphi(z) = \varphi(y) \in \varphi(Tz) ,$$

for some $y \in Tz \cap S(z)$. $\square$

3.4. **The equivalence of principles and completeness.** We prove first the equivalence between Ekeland, Takahashi and Caristi principles.

**Theorem 3.8.** *Let $(X, d, \preccurlyeq)$ be a preordered quasi-metric space and $\varphi : X \to \mathbb{R} \cup \{\infty\}$ a proper bounded below function. Then the following statements are equivalent.*

(wEk)  *The following holds*

(3.19)
$$\exists z \in \operatorname{dom} \varphi, \ \forall y \in S(z), \ \varphi(y) = \varphi(z).$$

(Tak)  *The following holds*

(3.20)
$$\begin{aligned}\{\forall x \in \operatorname{dom} \varphi, \ [\inf \varphi(X) < \varphi(x) \ \Rightarrow \ \exists y \in S(x), \ \varphi(y) < \varphi(x)]\} \\ \Longrightarrow \ \exists z \in \operatorname{dom} \varphi, \ \varphi(z) = \inf \varphi(X).\end{aligned}$$

(Car)  *If the mapping $T : X \to X$ satisfies*

(3.21)
$$Tx \in S(x) \ \text{for all} \ x \in \operatorname{dom} \varphi,$$

*then there exists $z \in \operatorname{dom} \varphi$ such that $\varphi(Tz) = \varphi(z)$ and $Tz \in \overline{\{z\}}$.*

*Proof.* Observe that, by Proposition 2.9, $y \in S(z)$ and $\varphi(y) = \varphi(z)$ implies $y \in \overline{\{z\}}$.

(wEk) $\iff$ (Tak).

We prove the equivalent assertion: $\neg$(wEk) $\iff \neg$ (Tak).

For convenience, denote by (Ta1) the expression
$$\forall x \in \operatorname{dom} \varphi, \ [\inf \varphi(X) < \varphi(x) \ \Rightarrow \ \exists y \in S(x), \ \varphi(y) < \varphi(x)].$$

The following equivalences hold.
$$\begin{aligned}\neg(\text{Tak}) \iff & (\text{Ta1}) \land (\forall z \in \operatorname{dom} \varphi, \ \varphi(z) > \inf \varphi(X)) \\ \iff & \forall x \in \operatorname{dom} \varphi, \ \exists y \in S(x), \ \varphi(y) < \varphi(x).\end{aligned}$$

By Proposition 2.9.(ii), $\varphi(y) \leq \varphi(x)$ for every $y \in S(x)$, so that
$$((y \in S(x)) \land (\varphi(y) \neq \varphi(x))) \iff ((y \in S(x)) \land (\varphi(y) < \varphi(x))).$$

But then

(3.22)
$$\begin{aligned}\neg(\text{wEk}) \iff & \forall x \in \operatorname{dom} \varphi, \ \exists y \in S(x), \ \varphi(y) \neq \varphi(x) \\ \iff & \forall x \in \operatorname{dom} \varphi, \ \exists y \in S(x), \ \varphi(y) < \varphi(x) \\ \iff & \neg(\text{Tak}).\end{aligned}$$

(wEk) $\Rightarrow$ (Car).

Suppose that $T : X \to X$ satisfies (3.21). By (wEk) there exists $z \in \operatorname{dom} \varphi$ such that $\varphi(x) = \varphi(z)$ for all $x \in S(z)$. Since, by hypothesis, $Tz \in S(z)$, it follows $\varphi(Tz) = \varphi(z)$ and $Tz \in \overline{\{z\}}$, because $S(z) \subseteq \overline{\{z\}}$.

$\neg$(wEk) $\Rightarrow \neg$(Car).

By (3.22), $\neg$(wEk) is equivalent to

(3.23)
$$\forall x \in \operatorname{dom} \varphi, \ \exists y_x \in S(x), \ \varphi(y_x) < \varphi(x).$$

Let $x_0 \in \operatorname{dom} \varphi$. Define $T : X \to X$ by $Tx = y_x$ for $x \in \operatorname{dom} \varphi$, where $y_x$ is given by (3.23), and $Tx = x_0$ for $x \in X \smallsetminus \{x_0\}$. Then $Tx \in S(x)$ but $\varphi(Tx) < \varphi(x)$, for all $x \in \operatorname{dom} \varphi$, i.e., (Car) fails. $\square$

Finally, we show that the validity of each of these principles is further equivalent to the sequential increasingly right $K$-completeness of the preordered quasi-metric space $X$. Various situations when the validity of a fixed point result or of a variational principle implies completeness are discussed in [13].

**Theorem 3.9.** *For a preordered quasi-metric space $(X, d, \preccurlyeq)$ with the order $\preccurlyeq$ satisfying (2.6) the following statements are equivalent.*

1. *The space $(X, d, \preccurlyeq)$ is sequentially increasingly right $K$-complete.*
2. (Ekeland variational principle - weak form) *For every proper bounded below increasingly lsc function $\varphi : X \to \mathbb{R} \cup \{\infty\}$ there exists $z \in \operatorname{dom} \varphi$ such that*

$$\varphi(x) = \varphi(z) \quad \text{for all} \ \ x \in S(z). \tag{3.24}$$

3. (Takahashi minimization principle) *Every proper bounded below increasingly lsc function $\varphi : X \to \mathbb{R} \cup \{\infty\}$ such that*

$$\varphi(x) > \inf \varphi(X) \ \Rightarrow \ \exists y \in S(x), \ \varphi(y) < \varphi(x),$$

*for every $x \in \operatorname{dom} \varphi$, attains its minimum on $X$.*

4. (Caristi fixed point theorem) *For every proper bounded below increasingly lsc function $\varphi : X \to \mathbb{R} \cup \{\infty\}$ and every mapping $T : X \to X$ such that*

$$Tx \in S(x) \ \ \text{for all} \ \ x \in \operatorname{dom} \varphi,$$

*there exists $z \in \operatorname{dom} \varphi$ such that $\varphi(Tz) = \varphi(z)$ and $Tz \in \overline{\{z\}}$.*

*Proof.* The equivalences $2 \iff 3 \iff 4$ are contained in Theorem 3.8 (even in a stronger form - with the same function $\varphi$)

The implication $1 \Rightarrow 2$ is contained in Theorem 3.3.

$2 \Rightarrow 1$.

The proof is inspired from [21]. We proceed by contradiction. Suppose that there exists an increasing right $K$-Cauchy sequence $(x_n)_{n \in \mathbb{N}}$ in $X$ which does not converge. Passing to a subsequence if necessary, we can further suppose that

$$d(x_{n+1}, x_n) < \frac{1}{2^{n+1}}, \quad \text{for all } n \in \mathbb{N}. \tag{3.25}$$

By Proposition 2.4), the sequence $(x_n)$ has no convergent subsequences, so we can also suppose that its terms are pairwise distinct. The same proposition implies that the set

$$B := \{x_n : n \in \mathbb{N}\}$$

is closed. Indeed, if it would exist $x \in \overline{B} \smallsetminus B$, then $x$ will be the limit of a subsequence of $(x_n)$, a contradiction.

Define $\varphi : X \to \mathbb{R}$ by

$$\varphi(x) = \begin{cases} \frac{1}{2^{n-1}} & \text{if} \ \ x = x_n \ \text{for some} \ \ n \in \mathbb{N}, \\ \infty & \text{for} \ \ x \in X \setminus B. \end{cases}$$

It follows that $\operatorname{dom} \varphi = B$ and that $\varphi$ is lsc (and so, *a fortiori*, increasingly lsc).

Indeed, the lsc of $\varphi$ is equivalent to the closedness of the set $[\varphi \leq b] := \{x \in X : \varphi(x) \leq b\}$, for every $b \in \mathbb{R}$. Since $[\varphi \leq b] = B$ for $b \geq 1$ and $[\varphi \leq b] = \emptyset$ for $b \leq 0$, it remains to examine the case $0 < b < 1$. If $k \in \mathbb{N}$ is such that
$$\frac{1}{2^k} \leq b < \frac{1}{2^{k-1}},$$
then, $[\varphi \leq b] = \{x_{k+1}, x_{k+2}, \dots\}$, and so it is closed.

If $x \in \operatorname{dom}\varphi = B$, then $x = x_k$, for some $k \in \mathbb{N}$. By (3.25),
$$d(x_{k+i}, x_k) \leq d(x_{k+i}, x_{k+i-1}) + d(x_{k+i-1}, x_{k+i-2}) + \cdots + d(x_{k+1}, x_k)$$
$$< \frac{1}{2^{k+i}} + \frac{1}{2^{k+i-1}} + \cdots + \frac{1}{2^{k+1}} < \frac{1}{2^k}.$$

Hence,
$$\varphi(x_{k+i}) + d(x_{k+i}, x_k) < \frac{1}{2^{k+i-1}} + \frac{1}{2^k} < \frac{1}{2^{k-1}} = \varphi(x_k),$$
that is, $x_{k+i} \in S(x_k)$ for all $i \in \mathbb{N}$. Since
$$\varphi(x_{k-i}) = \frac{1}{2^{k-i-1}} > \frac{1}{2^{k-1}} = \varphi(x_k),$$
for $i = 1, \dots, k-1$, it follows that $S(x_k) = \{x_k, x_{k+1}, x_{k+2}, \dots\}$.

The inequality $\varphi(x_{k+1}) = \frac{1}{2^k} < \frac{1}{2^{k-1}} = \varphi(x_k)$ shows that (wEk) fails. □

3.5. **The case of $T_1$ quasi-metric spaces.** Let us notice that a topological space $(X, \tau)$ is $T_1$ if and only if $\overline{\{x\}} = \{x\}$ for all $x \in X$. A quasi-metric space $(X, d)$ is $T_1$ (i.e., the topology $\tau_d$ is $T_1$) if and only if $d(x, y) > 0$ for all $x \neq y$ in $X$.

Taking into account these remarks, the results proved for arbitrary (i.e., $T_0$) quasi-metric spaces take the following forms in the $T_1$ case.

**Theorem 3.10.** *Suppose that $(X, d, \preccurlyeq)$ is a $T_1$ sequentially increasingly right $K$-complete preordered quasi-metric space such that the preorder $\preccurlyeq$ satisfies (2.6). Let $\varphi : X \to \mathbb{R} \cup \{\infty\}$ be a proper, bounded below and increasingly lsc function. The following are true.*

1. *(Ekeland variational principle - weak form, [10]) There exists $z \in \operatorname{dom}\varphi$ such that*

(3.26) $\qquad \varphi(z) < \varphi(x) + d(x, z) \quad$ *for all* $x \in \operatorname{dom}\varphi \setminus \{z\}$ *with* $z \preccurlyeq x$.

2. *(Takahashi principle, [2]) If for every $x \in \operatorname{dom}\varphi$,*

(3.27) $\qquad\qquad \varphi(x) > \inf \varphi(X) \;\Rightarrow\; S(x) \setminus \{x\} \neq \emptyset,$

   *then there exists $z \in \operatorname{dom}\varphi$ such that $\varphi(z) = \inf \varphi(X)$.*

3. *(Caristi fixed point theorem)*
    (a) *If the mapping $T : X \to X$ satisfies*

(3.28) $\qquad\qquad x \preccurlyeq Tx \;\; \text{and} \;\; d(Tx, x) + \varphi(Tx) \leq \varphi(x),$

   *for all $x \in \operatorname{dom}\varphi$, then there exists $z \in \operatorname{dom}\varphi$ such that $Tz = z$.*
    (b) *If $T : X \rightrightarrows X$ is a set-valued mapping such that*

(3.29) $\qquad\qquad S(x) \cap Tx \neq \emptyset,$

   *for every $x \in \operatorname{dom}\varphi$, then there exists $z \in \operatorname{dom}\varphi$ such that $z \in Tz$.*

*Proof.* 1. By Theorem 3.3 there exists $z \in \operatorname{dom}\varphi$ satisfying (3.3). Since $X$ is $T_1$, $\overline{\{z\}} = \{z\}$, so that $S(z) = \{z\}$. Taking into account this equality, (3.26) follows from (3.3).(ii) (with $y = z$).

2. If $y \in S(x) \smallsetminus \{x\}$, then $d(y,x) > 0$, so that
$$\varphi(y) < \varphi(y) + d(y,x) \leq \varphi(x),$$
This shows that condition (3.8) is verified by $\varphi$.

3. By Theorem 3.7, there exists $z \in \operatorname{dom}\varphi$ such that $\emptyset \neq Tz \cap \overline{\{z\}} = Tz \cap \{z\}$, which means that $z \in Tz$.

As we have noticed in the proof of the same theorem, the single-valued case follows from the multivalued one. $\square$

## References

[1] A. Alam and M. Imdad, *Relation-theoretic contraction principle*, J. Fixed Point Theory Appl. **17** (2015), no. 4, 693–702. 2

[2] S. Al-Homidan, Q. H. Ansari and G. Kassay, *Takahashi's minimization theorem and some related results in quasi-metric spaces*, J. Fixed Point Theory Appl. **21** (2019), no. 1, article 38, 20 p. 2, 17

[3] T. Q. Bao, S. Cobzaş, A. Soubeyran, *Variational principles, completeness and the existence of traps in behavioral sciences*, Ann. Oper. Res. **269** (2018), no. 1-2, 53–79. 2, 3, 5

[4] T. Q. Bao, B. S. Mordukhovich, and A. Soubeyran, *Fixed points and variational principles with applications to capability theory of wellbeing via variational rationality*, Set-Valued Var. Anal. **23** (2015), no. 2, 375–398. 2

[5] ______, *Minimal points, variational principles, and variable preferences in set optimization*, J. Nonlinear Convex Anal. **16** (2015), no. 8, 1511–1537. 2

[6] ______, *Variational analysis in psychological modeling*, J. Optim. Theory Appl. **164** (2015), no. 1, 290–315. 2

[7] T. Q. Bao and A. Soubeyran, *Variational analysis in cone pseudo-quasimetric spaces and applications to group dynamics*, J. Optim. Theory Appl. **170** (2016), no. 2, 458–475. 2

[8] T. Q. Bao and M. A. Théra, *On extended versions of Dancs-Hegedüs-Medvegyev's fixed-point theorem*, Optimization **66** (2017), 875–887. 3, 7

[9] G. C. Bento, J. X. Cruz Neto, J. O. Lopes, P. A. Soares, Jr., and A. Soubeyran, *Generalized proximal distances for bilevel equilibrium problems*, SIAM J. Optim. **26** (2016), no. 1, 810–830. 2

[10] S. Cobzaş, *Completeness in quasi-metric spaces and Ekeland Variational Principle*, Topology Appl., **158** (2011), pp. 1073-1084. 2, 17

[11] ______, *Ekeland variational principle in asymmetric locally convex spaces*, Topology Appl. **159** (2012), no. 10-11, 2558–2569. 2

[12] ______, *Functional Analysis in Asymmetric Normed Spaces*, Frontiers in Mathematics, Birkhäuser/Springer Basel AG, Basel, 2013. 3, 4

[13] ______, *Fixed points and completeness in metric and in generalized metric spaces*, Fundam. Prikl. Mat. **22** (2018), no. 1, 127–215 (Russian), (see also: arXiv:1508.05173v4, 71 p). 16

[14] ______, *Ekeland, Takahashi and Caristi principles in quasi-pseudometric spaces*, Topology Appl. **265** (2019), 106831, 22. 2, 6, 11

[15] S. Dancs, M. Hegedüs, and P. Medvegyev, *A general ordering and fixed-point principle in complete metric space*, Acta Sci. Math. (Szeged), 46 (1983), pp. 381-388. 3, 7

[16] I. Ekeland, *Sur les problèmes variationnels*, C. R. Acad. Sci. Paris Sér. A-B, 275 (1972), pp. 1057-1059. 1

[17] ______, *On the variational principle*, J. Math. Anal. Appl. 47 (1974), 324–353. 1

[18] ______, *Nonconvex minimization problems*, Bull. Amer. Math. Soc. (N.S.) 1 (1979), no. 3, 443–474. 1

[19] A. Hamel, *Remarks to an equivalent formulation of Ekeland's variational principle*, Optimization **31** (1994), 233–238. 1
[20] J. Jachymski, *The contraction principle for mappings on a metric space with a graph*, Proc. Amer. Math. Soc. **136** (2008), no. 4, 1359–1373. MR 2367109 2, 5
[21] E. Karapınar and S. Romaguera, *On the weak form of Ekeland's variational principle in quasi-metric spaces*, Topology Appl. **184** (2015), 54–60. 2, 6, 16
[22] J. C. Kelly, *Bitopological spaces*, Proc. London Math. Soc. **13** (1963), 71–89. 3
[23] H.-P. A. Künzi, M. Mršević, I. L. Reilly, and M. K. Vamanamurthy, *Convergence, precompactness and symmetry in quasi-uniform spaces*, Math. Japon. **38** (1993), 239–253. 5
[24] I. Meghea, *Ekeland variational principle. With generalizations and variants*, Old City Publishing, Philadelphia, PA; Éditions des Archives Contemporaines, Paris, 2009. 2
[25] J. J. Nieto, R. L. Pouso, and R. Rodríguez-López, *Fixed point theorems in ordered abstract spaces*, Proc. Amer. Math. Soc. **135** (2007), no. 8, 2505–2517. MR 2302571 2
[26] J. J. Nieto and R. Rodríguez-López, *Contractive mapping theorems in partially ordered sets and applications to ordinary differential equations*, Order **22** (2005), no. 3, 223–239 (2006). MR 2212687 2, 5
[27] ______, *Existence and uniqueness of fixed point in partially ordered sets and applications to ordinary differential equations*, Acta Math. Sin. (Engl. Ser.) **23** (2007), no. 12, 2205–2212. MR 2357454 2
[28] J.-P. Penot, *Fixed point theorems without convexity*, Bull. Soc. Math. Fr., Suppl., Mém. **60** (1979), 129-152 . 6
[29] ______, *The drop theorem, the petal theorem and Ekeland's variational principle*, Nonlinear Anal. **10** (1986), no. 9, 813–822. 1
[30] A. Petruşel and I. A. Rus, *Fixed point theorems in ordered L-spaces*, Proc. Amer. Math. Soc. **134** (2006), no. 2, 411–418. MR 2176009 2
[31] A. Ran and M. C. B. Reurings, *A fixed point theorem in partially ordered sets and some applications to matrix equations*, Proc. Amer. Math. Soc. **132** (2004), no. 5, 1435–1443. MR 2053350 2
[32] J.-H. Qiu, F. He, and A. Soubeyran, *Equilibrium versions of variational principles in quasi-metric spaces and the robust trap problem*, Optimization **67** (2018), no. 1, 25–53. 2
[33] I. L. Reilly, P. V. Subrahmanyam, M. K. Vamanamurthy, *Cauchy sequences in quasi-pseudo-metric spaces,* Monatsh. Math. **93** (1982), 127–140. 4
[34] S. Shukla and R. Rodríguez-López, *Fixed points of multi-valued relation-theoretic contractions in metric spaces and application*, Quaest. Math. **43** (2020), no. 3, 409–424. 2
[35] W. Takahashi, *Existence theorems generalizing fixed point theorems for multivalued mappings*, Fixed point theory and applications (Marseille, 1989), Pitman Res. Notes Math. Ser., vol. 252, Longman Sci. Tech., Harlow, 1991, pp. 397–406. 1
[36] ______, *Nonlinear functional analysis. Fixed point theory and its applications*, Yokohama Publishers, Yokohama, 2000. 1

Babeş-Bolyai University, Faculty of Mathematics and Computer Science, 400 084 Cluj-Napoca, Romania

*Email address*: scobzas@math.ubbcluj.ro